\documentclass[12pt]{amsart}
\renewcommand{\bar}{\overline}
\usepackage{amsmath,amsthm,amscd,euscript,amssymb}
\usepackage{graphicx}
\usepackage{epstopdf}
\usepackage{color}
\usepackage[most]{tcolorbox}

\newcommand{\myexample}[2]{
   \begin{tcolorbox}[breakable,colback=black!5!white,colframe=black,title={%Example
   #1}]
        #2
    \end{tcolorbox}
}

\renewcommand{\bar}{\overline}
\def \r{\mathbb R}
\def \c{\mathbb C}
\def \q{\mathbb Q}
\def \z{\mathbb Z}

\def \M{\mathcal M}

\def \GL{\hbox{\rm GL}}
\def \SL{\hbox{\rm SL}}

\def \re{\hbox{\rm Re}}
\def \im{\hbox{\rm Im}}

\newtheorem{theorem}{Theorem}[section]
\newtheorem{lemma}[theorem]{Lemma}

\newtheorem{proposition}[theorem]{Proposition}
\newtheorem{corollary}[theorem]{Corollary}
\theoremstyle{remark}
\newtheorem{remark}[theorem]{Remark}
\theoremstyle{definition}
\newtheorem{definition}[theorem]{Definition}
\newtheorem{example}[theorem]{Example}
\newtheorem{problem}{Problem}
\newtheorem{conjecture}[problem]{Conjecture}

\title[The $\sin^2$-algorithm in the non-totally-real case]{On periodicity of the $\sin^2$-algorithm for non-totally-real cubic vectors}

\author{Oleg Karpenkov}

\date{25 November 2026}


\email[Oleg Karpenkov]{karpenk@liverpool.ac.uk}

\begin{document}
\IfFileExists{epsf.tex}{\input{epsf}}{}

\begin{abstract}
The $\sin^2$-algorithm is a modification of the Jacobi-Perron algorithm in the three dimensional case;
its periodicity is known for totally-real conjugate cubic vectors.
In this paper we extend the $\sin^2$-algorithm to the case of cubic vectors with complex conjugate vectors:
the angle between two real planes is replaced by the Hermitian angle between two complex planes.
We prove that the extended algorithm is eventually periodic for all non-totally-real cubic vectors.
This provides an answer in the non-totally-real case to the question of algebraic periodicity for cubic irrationalities posed in 1848 by Ch.~Hermite.
One of the statements used in the proof is a geometric inequality whose proof is computer-assisted:
it is verified on 16507 boxes in exact rational arithmetic.
This inequality is also formalised in the proof assistant Lean.
We show that it implies the similar inequality in the totally-real case,
which gives another proof of the periodicity of the original $\sin^2$-algorithm.
\end{abstract}

\maketitle
\tableofcontents

\section*{Introduction}

Recall that a cubic number is a root of a cubic polynomial with integer coefficients
irreducible over the field $\q$.
A cubic number is called {\it totally-real} if the cubic integer polynomial defining it
has three real roots; otherwise it is {\it non-totally-real}.
We say that a vector in $\r^3$ is {\it cubic} if its coordinates span a cubic field.

\vspace{2mm}

In 1770 J.-L.~Lagrange showed in~\cite{Lagrange1770} that
a real number is a root of an irreducible quadratic polynomial
if and only if its continued fraction is periodic.
The question of similar periodic representations of cubic irrationalities was
posed by Ch.~Hermite in 1848 (see e.g.~\cite{Picard1901}, \cite{Hermite1850}).
The first generalization of the Euclidean algorithm to the three dimensional case
was introduced by C.G.J.~Jacobi in~\cite{Jacobi1868} and further developed by O.~Perron in~\cite{Perron1907}.
The Jacobi-Perron algorithm is periodic for some cubic vectors, however it is believed
to be non-periodic for some others (see, e.g.,~\cite{Schweiger2000} and Chapter~23.4 in~\cite{KarpenkovGCF2013}).

\vspace{2mm}

In~\cite{Karpenkov2021} we have introduced the $\sin^2$-algorithm.
It works with a triple of vectors: a vector $\xi$ and its two conjugate vectors $\nu_1$ and $\nu_2$.
At each step it picks the Jacobi-Perron type transformation that maximizes the sine of the angle between the planes spanned by
$(\xi,\nu_1)$ and by $(\xi,\nu_2)$ in the coordinates of the new basis.
The main result of~\cite{Karpenkov2021} is the periodicity of the $\sin^2$-algorithm for totally-real cubic vectors.
The non-totally-real case was left open (see the last section of~\cite{Karpenkov2021} and the discussion in~\cite{Karpenkov2021-2}),
with a suggestion to use the analytic extension of the formula for $\sin^2\alpha$ to complex numbers.

\vspace{2mm}

The aim of the present paper is to show that this suggestion works.
For a real vector $\xi$ and a pair of complex conjugate vectors $\nu$, $\bar\nu$ we consider
the {\it Hermitian} angle $\alpha$ between the complex planes spanned by $(\xi,\nu)$ and by $(\xi,\bar\nu)$.
The value $\sin^2\alpha$ is given by a simple rational expression
(the analytic extension of the real formula equals $-\tan^2\alpha$, see Remark~\ref{analytic}).
Then the algorithm is literally the same as in~\cite{Karpenkov2021};
we call it the {\it complex $\sin^2$-algorithm}.
The main result of this paper is as follows (see Theorems~\ref{eventually-periodic-c} and~\ref{eventually-periodic-algorithm-c}).

\vspace{2mm}

{\noindent
{\bf Main result.}
{\it
The complex $\sin^2$-algorithm is eventually periodic for a triple of vectors $(\xi,\nu,\bar\nu)$ in $\z$-general position
if and only if $(\xi,\nu,\bar\nu)$ is a triple of conjugate cubic vectors.
}}

\vspace{2mm}

The proof follows the scheme of the totally-real case, however its two ingredients are obtained in a different way.

\vspace{1mm}

{\noindent
{\it $($i$)$ Finiteness of $\varepsilon$-caps} (Theorem~\ref{eventually-periodic-c}$($i$)$).}
For every $\varepsilon>0$ there are only finitely many separating states of a given cubic vector
with $\sin^2\alpha\ge\varepsilon$.
This is a statement of geometry of numbers; its proof is short and uses only the
norm form and the action of the Dirichlet group.
It is formalised in the proof assistant Lean~4 as well (see Remark~\ref{lean-caps}).

\vspace{1mm}

{\noindent
{\it $($ii$)$ Absence of flat local maxima} (Theorem~\ref{key}).}
If $\sin^2\alpha<\varepsilon_0=400/10201$ for some state (cubic or not), then there exists an admissible JP-transformation that increases $\sin^2\alpha$.
This is a statement of elementary geometry:
it concerns a rational function of four real variables and a discrete family of transformations, no arithmetic is involved.
We give a computer-assisted proof.
After a compactification of the space of parameters, the inequality is verified on $16507$ boxes by interval arithmetic
over $\q$ (so there are no rounding errors).
In addition, this statement and its proof are formalised in the proof assistant Lean~4 (see Remark~\ref{lean}).

\vspace{2mm}

These two statements imply that the set of states where the algorithm decreases $\sin^2\alpha$ is finite,
and the periodicity follows as in~\cite{Karpenkov2021}.

\vspace{2mm}

Finally, we show that Theorem~\ref{key} implies the same inequality for the angle between two real planes
(Theorem~\ref{key-real}).
So in the totally-real case the finiteness of the extremal sets, which is the most complicated part of~\cite{Karpenkov2021},
follows from the finiteness of $\varepsilon$-caps (see Corollary~\ref{real-omega}).

\vspace{2mm}

{\noindent
{\bf This paper is organized as follows.}
}
In Section~\ref{sec-def} we define the complex $\sin^2$-algorithm, formulate the main results, and show several examples.
In Section~\ref{sec-caps} we prove the finiteness of $\varepsilon$-caps.
Section~\ref{sec-key} is dedicated to the proof of Theorem~\ref{key} and to its totally-real analogue.
In Section~\ref{sec-proofs} we deduce the periodicity of the algorithm.
We conclude the paper in Section~\ref{sec-open} with several open questions.

\section{Definition of the complex $\sin^2$-algorithm and theorems on its periodicity}
\label{sec-def}

We follow the definitions and the notation of~\cite{Karpenkov2021} as close as possible.

\subsection{Non-totally-real cubic states}

\begin{definition}\label{state-c}
A {\it $($complex$)$ state} is a collection of coordinates for three vectors $(\xi,\nu,\bar\nu)$
linearly independent over $\c$, where $\xi=(x,y,z)$ is a real vector satisfying
$$
x\ge y\ge z>0,
$$
and $\nu\in\c^3$.
We consider states up to a multiplication of $\xi$ by a positive real number
and of $\nu$ by a non-zero complex number.
\end{definition}

\begin{definition}
We say that a state $s=(\xi,\nu,\bar\nu)$ is a {\it non-totally-real cubic} state if
there exists a matrix  $B\in\SL(3,\z)$ whose characteristic polynomial is irreducible over $\q$
and has exactly one real root, such that the vectors $\xi,\nu,\bar\nu$ are eigenvectors of $B$.
\end{definition}

As in the totally-real case, every real cubic vector $\xi$ whose conjugate vectors are not real is an eigenvector
of such a matrix (this follows from Dirichlet's unit theorem, see~\cite{Borevich1966} and Chapter~17 in~\cite{KarpenkovGCF2013}).
The conjugate vectors $\nu$ and $\bar\nu$ are then determined by $\xi$ up to multiplication by complex numbers.

\subsection{Hermitian angle between two complex planes}

For two vectors $u=(u_1,u_2,u_3)$ and $v=(v_1,v_2,v_3)$ in $\c^3$
we write $u\times v$ for their vector product (defined by the usual formula)
and
$$
\langle u,v\rangle=u_1\bar v_1+u_2\bar v_2+u_3\bar v_3, \qquad |u|^2=\langle u,u\rangle.
$$
Recall that the {\it Hermitian angle} between two complex planes in $\c^3$ that are
orthogonal to vectors $m_1$ and $m_2$ with respect to the Hermitian product
is the angle $\alpha\in[0,\pi/2]$ defined by
$$
\cos\alpha=\frac{|\langle m_1,m_2\rangle|}{|m_1|\cdot|m_2|}.
$$

\begin{definition}\label{sin-def}
Consider a state $s=(\xi,\nu,\bar\nu)$.
Denote by $\alpha(s)$ the Hermitian angle between the complex planes spanned by
pairs of vectors
$(\xi,\nu)$ and $(\xi,\bar\nu)$ respectively.
\end{definition}

The plane spanned by $(\xi,\nu)$ is Hermitian-orthogonal to the vector $\bar n$, where $n=\xi\times \nu$,
and the plane spanned by $(\xi,\bar\nu)$ is Hermitian-orthogonal to $n$.
Therefore, for $n=(n_1,n_2,n_3)$ we have
\begin{equation}\label{sin-n}
\sin^2\alpha(s)=1-\frac{|n_1^2+n_2^2+n_3^2|^2}{\big(|n_1|^2+|n_2|^2+|n_3|^2\big)^2}.
\end{equation}
This value does not change if we multiply $\xi$ and $\nu$ by non-zero numbers or swap $\nu$ and $\bar \nu$.

\begin{remark}\label{analytic}
For three real vectors $(\xi,\nu_1,\nu_2)$ the angle $\alpha$ between the planes spanned by $(\xi,\nu_1)$ and $(\xi,\nu_2)$
satisfies
$$
\sin^2\alpha=1-\frac{(n\cdot n')^2}{(n\cdot n)(n'\cdot n')}, \quad \hbox{where $n=\xi\times\nu_1$, $n'=\xi\times\nu_2$,}
$$
and the dot denotes the bilinear scalar product.
If we substitute here $\nu_1=\nu$ and $\nu_2=\bar\nu$, we get
$$
1-\frac{\big(|n_1|^2+|n_2|^2+|n_3|^2\big)^2}{|n_1^2+n_2^2+n_3^2|^2}=-\tan^2\alpha(s).
$$
So the analytic extension of the formula for $\sin^2\alpha$ suggested in~\cite{Karpenkov2021} is
a decreasing function of $\sin^2\alpha(s)$:
its minimization is equivalent to the maximization of the Hermitian $\sin^2\alpha(s)$.
\end{remark}

It is convenient to write $\sin^2\alpha(s)$ in terms of the dual row.

\begin{proposition}\label{sin-ell}
Let $s=(\xi,\nu,\bar\nu)$ be a state and let $\ell=(\ell_1,\ell_2,\ell_3)$ be a non-zero complex row
satisfying $\ell\xi=0$ and $\ell\bar\nu=0$ $($here vectors are columns$)$. Write $\ell=a+ib$ with real rows $a$ and $b$.
Then
\begin{equation}\label{sin-l}
\sin^2\alpha(s)=1-\frac{|\ell_1^2+\ell_2^2+\ell_3^2|^2}{\big(|\ell_1|^2+|\ell_2|^2+|\ell_3|^2\big)^2}
=\frac{4|a\times b|^2}{\big(|a|^2+|b|^2\big)^2}.
\end{equation}
In particular, $0<\sin^2\alpha(s)\le 1$.
\end{proposition}

\begin{proof}
The row $\ell$ is defined up to a non-zero complex factor, and it is proportional to $\bar n$.
This gives the first equality.
Further
$$
\big(|a|^2+|b|^2\big)^2-\big|\,|a|^2-|b|^2+2i\, a\cdot b\,\big|^2=4\big(|a|^2|b|^2-(a\cdot b)^2\big)=4|a\times b|^2.
$$
Finally, the rows $\ell$ and $\bar\ell$ are linearly independent, since $\bar \ell\nu=0$ while $\ell\nu\ne 0$.
Hence $a$ and $b$ are linearly independent and $a\times b\ne 0$.
\end{proof}

\begin{remark}
Consider the linear map $\r^3\to\r^2$ whose matrix has the rows $a$ and $b$;
it sends the three coordinate vectors to the planar vectors $\ell_1,\ell_2,\ell_3\in\c=\r^2$.
The image of the unit ball under this map is bounded by an ellipse.
Let $\sigma_1\ge\sigma_2>0$ be the semi-axes of this ellipse, and let $t=\sigma_2/\sigma_1$.
Then
$$
\sigma_1^2+\sigma_2^2=|a|^2+|b|^2 \qquad \hbox{and} \qquad \sigma_1\sigma_2=|a\times b|,
$$
and hence by~\eqref{sin-l} we have
$$
\sin^2\alpha(s)=\frac{4t^2}{(1+t^2)^2}.
$$
So $\sin^2\alpha(s)=1$ if and only if the ellipse is a circle, which is equivalent to $\ell_1^2+\ell_2^2+\ell_3^2=0$
$($for instance, this is the case for three vectors of the same length with the angles $2\pi/3$ between them$)$.
The value $\sin^2\alpha(s)$ is small if and only if the ellipse is thin, i.e.,
the three vectors $\ell_1$, $\ell_2$, and $\ell_3$ are almost collinear.
In the latter case we say that the state is {\it almost flat}.
\end{remark}

\subsection{Supporting and separating bases}

A vector is called {\it integer} if all its coordinates are integers.
A basis is {\it integer} if it generates the lattice of integer points.
Recall that a basis $E$ is {\it supporting} for $\xi$
if the non-negative orthant of $E$ contains either $\xi$ or $-\xi$.

For a triple $(\xi,\nu,\bar\nu)$ as above denote by $\pi(\nu)$ the real plane spanned by $\re\,\nu$ and $\im\,\nu$.
This plane does not contain $\xi$.

\begin{definition}\label{separating-c}
A supporting basis $E$ for $\xi$ is said to be {\it separating} for $(\xi,\nu,\bar\nu)$
if the non-negative orthant of $E$ meets the plane $\pi(\nu)$ only at the origin.
We say that a state $s$ is {\it separating} if the coordinate basis is separating for the vectors of $s$.
\end{definition}

Denote by $L_1$ the real linear form that vanishes on $\pi(\nu)$ and satisfies $L_1(\xi)=1$.

\begin{proposition}\label{sep-heights}
Let $E=(e_1,e_2,e_3)$ be a basis whose non-negative orthant contains $\xi$.
Then $E$ is separating for $(\xi,\nu,\bar\nu)$ if and only if
$$
L_1(e_1)>0, \quad L_1(e_2)>0, \quad \hbox{and} \quad L_1(e_3)>0.
$$
\end{proposition}

\begin{proof}
If all three values are positive, then $L_1$ is positive at all non-zero points of the non-negative orthant.
Conversely, if the orthant meets $\pi(\nu)$ only at the origin, then $L_1$ does not vanish on the orthant without the origin, which is a connected set.
Since $L_1(\xi)>0$, the form $L_1$ is positive on it.
\end{proof}

\begin{remark}
In the totally-real case a supporting basis is separating if its non-negative orthant
does not include both $\pm\nu_1$ and $\pm\nu_2$.
Definition~\ref{separating-c} is the counterpart of this condition for the invariant plane $\pi(\nu)$, which replaces the two real eigenlines.
\end{remark}

\subsection{Admissible JP-transformations}

The transformations are the same as in~\cite{Karpenkov2021}.
Let $E=(e_1,e_2,e_3)$ be some basis.
\begin{itemize}
\item
Denote by
$V_{\alpha,\beta;\gamma}(E)$ the basis
$$
(e_1,e_2+\gamma e_1,e_3+\alpha e_1+\beta e_2),
$$
where $\alpha$, $\beta$, $\gamma$ are non-negative integers, not all equal to zero.

\item
Denote by
$W(E)$ the basis
$$
(e_1+e_3,e_2+e_1, e_3).
$$
\end{itemize}
These basis transformations are called {\it JP-transformations}.
A JP-transformation $M$ is  {\it admissible} for the pair $(\xi, E)$  if
$M(E)$ is a supporting basis for~$\xi$.
For a vector $\xi$ with a supporting basis $E$ we
denote by $\M_\xi$ the set of all admissible JP-transformations.
If the coordinates of $\xi$ in $E$ are $(x,y,z)$, then its coordinates in the new bases are
\begin{equation}\label{new-coord}
\begin{array}{l}
V_{\alpha,\beta;\gamma}:\quad \big(x-\alpha z-\gamma(y-\beta z),\, y-\beta z,\, z\big);\\
W:\quad \big(x - y,\, y,\, z-(x-y)\big).
\end{array}
\end{equation}

\begin{proposition}\label{JP-prop}
Let $\xi$ be a non-zero vector with a supporting basis $E$,
such that the coordinates $(x,y,z)$ of $\xi$
in the basis $E$ satisfy $x\ge y\ge z>0$.
Then the following statements hold.
\\
{\it $($i$)$} The JP-transformation $V_{\alpha,\beta;\gamma}$ is admissible if and only if
$$
y-\beta z\ge 0 \quad \hbox{and} \quad x-\alpha z-\gamma(y-\beta z)\ge 0.
$$
{\it $($ii$)$} The JP-transformation $W$  is admissible if and only if $z\ge x-y$.
\\
{\it $($iii$)$} For every JP-transformation $M$ the non-negative orthant of $M(E)$ is contained in the non-negative orthant of $E$.
\\
{\it $($iv$)$} Let $E$ be a separating basis for $(\xi,\nu,\bar\nu)$
and let $M$ be an admissible JP-transformation for $\xi$.
Then $M(E)$ is a separating basis for $(\xi,\nu,\bar\nu)$.
\end{proposition}

\begin{proof}
The first two items follow from~\eqref{new-coord}.
The third item is straightforward: the vectors of $M(E)$ are non-negative integer combinations of the vectors of $E$.
The last item follows from the third one.
\end{proof}

\begin{remark}
If the coordinates of $\xi$ are linearly independent over $\q$ (which is the case for cubic vectors),
then all the inequalities of Proposition~\ref{JP-prop} are strict for admissible transformations
and we arrive at the conditions of~\cite{Karpenkov2021}.
In this case the set $\M_\xi$ is finite.
\end{remark}

\subsection{Definition of the map $\Phi$ and the complex $\sin^2$-algorithm}

For a state $s$ and a transformation $M\in\M_\xi$ we denote by $M(s)$
the collection of coordinates of the vectors of $s$ in the basis $M(E)$, where $E$ is the coordinate basis.

\begin{definition}\label{optimization rule-c}
Let $s$ be a state.
Let  $M\in \M_\xi$ be
the transformation with the greatest possible value of $\sin^2\alpha(M(s))$.
Let
$$
\Phi(s)=T\circ M,
$$
where $T$ is a transposition of the basis vectors that puts the coordinates of $\xi$
in decreasing order.
We call the transformation $\Phi(s)$ the {\it $\sin^2$-transformation} for $s$.
By $\Phi(s)$ we also denote the resulting state.
\end{definition}

\begin{remark}\label{ties}
If the maximum is attained at several transformations,
we take the first of them with respect to the lexicographic
ordering of the triples $(\alpha,\beta,\gamma)$ for the elements $V_{\alpha,\beta;\gamma}$,
additionally setting $W>V_{*,*;*}$ (as in~\cite{Karpenkov2021}).
In the non-totally-real cubic case this does happen: see Example~\ref{ex-tie} below.
\end{remark}

\myexample{Complex $\sin^2$-algorithm}
{
{\noindent
{\bf \underline{Input data}:} We are given a real vector $\xi$ and a complex vector $\nu$ such that $\xi$, $\nu$, and $\bar\nu$ are linearly independent in $\c^3$.
We start with the coordinate basis $E$.
}
\\
{\bf \underline{Preliminary Stage 1: Finding a supporting basis}.}
First we multiply the vectors of $E$ by $\pm 1$ and swap them such that the coordinates
of the vector $\xi$ in the new basis satisfy
$$
x\ge y\ge z>0.
$$

\vspace{2mm}

{\noindent
{\bf \underline{Preliminary Stage 2: Finding a separating basis}.}
While the current basis $(e_1,e_2,e_3)$ is not separating for $(\xi,\nu,\bar\nu)$, we replace it with the basis
$$
\Big(e_2+\Big\lfloor\frac{x}{y}\Big\rfloor e_1+\Big\lfloor\frac{z}{y}\Big\rfloor e_3,\; e_3,\; e_1\Big),
$$
where $(x,y,z)$ are the coordinates of $\xi$ in the current basis
$($this is a step of the Jacobi-Perron algorithm: the new coordinates of $\xi$ are
$(y,\,z-\lfloor z/y\rfloor y,\,x-\lfloor x/y\rfloor y)$; the vectors are not reordered between these steps$)$.
Once a separating basis is reached, we reorder its vectors as in Preliminary Stage~$1$.
}

\vspace{1mm}

{\noindent
{\bf \underline{Main Stage}:}
Each step we perform the following iteration: $s \to \Phi(s)$.
}

{\noindent
{\bf \underline{Termination of the algorithm}:}
Once we produce a triple $(x,y,z)$ with one of the coordinates equal to 0,
the algorithm terminates.
}

\vspace{1mm}

{\noindent
{\bf \underline{Output}:}
A sequence of transformations generated in the iterations of
Preliminary Stages~1 and~2 and Main Stage.
}
}

\vspace{2mm}

\begin{proposition}\label{prelim}
Let the coordinates of $\xi$ be linearly independent over $\q$.
Then Preliminary Stage~$2$ terminates and produces a separating basis for $(\xi,\nu,\bar\nu)$.
\end{proposition}

\begin{proof}
All the bases of Preliminary Stage~2 are integer bases whose non-negative orthants contain $\xi$.
Since the coordinates of $\xi$ are independent over $\q$, the Jacobi-Perron algorithm does not terminate
$($all the coordinates stay positive$)$.
By the theorem of O.~Perron on the convergence of the Jacobi-Perron algorithm (see~\cite{Perron1907} and~\cite{Schweiger2000})
the directions of all three vectors of the $n$-th basis tend to the direction of $\xi$ as $n$ tends to infinity.
Since $L_1(\xi)>0$, the form $L_1$ is positive at all three vectors of the $n$-th basis for sufficiently large $n$.
By Proposition~\ref{sep-heights} such a basis is separating.
\end{proof}

\begin{remark}
In all the examples that we have computed the main stage itself brings a supporting basis to a separating basis
in at most five steps, so Preliminary Stage~2 seems to be excessive
$($compare with the similar conjecture on Preliminary Stage~2 in~\cite{Karpenkov2021}$)$.
\end{remark}

\subsection{Periodicity of the complex $\sin^2$-algorithm}
\label{Periodicity-c}

\begin{definition}
Let $(\xi,\nu,\bar\nu)$ be a triple of vectors as above.
\begin{itemize}
\item
We say that a state $s$ is a {\it $\xi$-state} if
in some integer basis the coordinates of vectors $\xi$, $\nu$, and $\bar\nu$ are respectively proportional to the vectors of $s$.

\item
The set of all separating $\xi$-states $s$ where $\sin^2\alpha(s)\ge\varepsilon$ is said to be
the {\it $\varepsilon$-cap} for $(\xi,\nu,\bar\nu)$ and denoted by $\Omega_\varepsilon(\xi;\nu,\bar\nu)$.

\item
The set of all separating $\xi$-states $s$ satisfying
$$
\sin^2\alpha(s)\ge\sin^2\alpha(\Phi(s))
$$
is said to be
{\it extremal} for $\xi$ and denoted by $\Omega_{\max}(\xi;\nu,\bar\nu)$.
\end{itemize}
\end{definition}

The first main statement is a general inequality.
It holds for all states, no algebraicity is assumed.

\begin{theorem}\label{key}
Let
$$
\varepsilon_0=\frac{400}{10201}=0.0392\ldots
$$
Then for every state $s=(\xi,\nu,\bar\nu)$ satisfying $\sin^2\alpha(s)<\varepsilon_0$
there exists an admissible JP-transformation $M\in\M_\xi$ such that
$$
\sin^2\alpha(M(s))>\sin^2\alpha(s).
$$
\end{theorem}

We prove Theorem~\ref{key} in Section~\ref{sec-key}. The proof is computer-assisted.

\begin{theorem}\label{eventually-periodic-c}
Let $B\in \SL(3,\z)$ be a matrix whose characteristic polynomial is irreducible over $\q$
and has exactly one real root.
Let also $(\xi,\nu,\bar\nu)$ be some eigenbasis of $B$, where $\xi$ is real.
Then the following two statements hold.

{\noindent
{\it $($i$)$}}
For every $\varepsilon>0$ the set $\Omega_\varepsilon(\xi;\nu,\bar\nu)$ is finite.

{\noindent
{\it $($ii$)$}}
The set $\Omega_{\max}(\xi;\nu,\bar\nu)$ is finite.

\end{theorem}

\begin{proof}
Item~$($i$)$ is proven in Section~\ref{sec-caps}.
By Theorem~\ref{key} we have
$$
\Omega_{\max}(\xi;\nu,\bar\nu)\subset \Omega_{\varepsilon_0}(\xi;\nu,\bar\nu).
$$
Therefore, Item~$($ii$)$ follows from Item~$($i$)$.
\end{proof}

\begin{corollary}\label{period-c}
Let $B$ and $(\xi,\nu,\bar\nu)$ be as in Theorem~\ref{eventually-periodic-c}.
Suppose that the coordinate basis is separating, and let $s_0$ be the corresponding state.
Then the sequence of states $s_n=\Phi^n(s_0)$ is eventually periodic.
\end{corollary}

We say that the complex $\sin^2$-algorithm is {\it eventually periodic} for $(\xi,\nu,\bar\nu)$
if it does not terminate and the sequence of states of its main stage is eventually periodic.
In this case the sequence of transformations in the output is eventually periodic as well.
We say that the vectors $\xi$, $\nu$, and $\bar\nu$ are in {\it $\z$-general position}
if any two of them span a complex plane that does not have non-zero integer vectors in it.
Triples of cubic conjugate vectors are always in $\z$-general position.

\begin{theorem}\label{eventually-periodic-algorithm-c}
Consider a real vector $\xi$ and a complex vector $\nu$ such that the triple $(\xi,\nu,\bar\nu)$ is linearly independent and in $\z$-general position.
Then the complex $\sin^2$-algorithm is eventually periodic for $(\xi,\nu,\bar\nu)$
if and only if $(\xi,\nu,\bar\nu)$ is a triple of conjugate cubic vectors.
\end{theorem}

Corollary~\ref{period-c} and Theorem~\ref{eventually-periodic-algorithm-c} are proven in Section~\ref{sec-proofs}.

\begin{remark}
In the totally-real case the finiteness of the extremal set $\Omega_{\max}$ is the most complicated part of~\cite{Karpenkov2021}.
In Subsection~\ref{sec-real} we show that Theorem~\ref{key} implies its totally-real analogue,
and hence the finiteness of $\Omega_{\max}$ in the totally-real case as well.
\end{remark}

\subsection{Examples}

In the examples below $\theta$ is the real root of a cubic polynomial, $\omega$ is one of its complex roots,
$\xi=(\theta^2,\theta,1)$, and $\nu=(\omega^2,\omega,1)$.
The coordinate basis is separating in all of them.

\begin{example}\label{ex-2}
Let $\theta=\sqrt[3]{2}$. The algorithm is purely periodic with period of length $3$:
$$
\renewcommand{\arraystretch}{1.2}
\begin{array}{|c|c|c|c|}
\hline
n & \hbox{coordinates of $\xi$, normalized} & \sin^2\alpha(s_n) & \hbox{JP-transformation}\\
\hline
0 & (1.587401,\ 1.259921,\ 1) & 0.932441 & V_{1,1;1}\\
\hline
1 & (3.847322,\ 1.259921,\ 1) & 0.835690 & V_{1,0;0}\\
\hline
2 & (2.847322,\ 1.259921,\ 1) & 0.938369 & V_{0,0;1}\\
\hline
\end{array}
$$
and $s_3=s_0$.
The basis after one period consists of the vectors $(1,1,1)$, $(2,1,1)$, and $(2,2,1)$.
The corresponding matrix has the characteristic polynomial $x^3-3x^2-3x-1$;
its eigenvalue at $\xi$ is the unit $1+\sqrt[3]{2}+\sqrt[3]{4}$.
\end{example}

\begin{example}
Let $\theta$ be the real root of $x^3-x-1$.
Here $s_0=(1.754878,\,1.324718,\,1)$, $\sin^2\alpha(s_0)=0.777999$, and
the $\sin^2$-transformation is $V_{1,0;0}$ followed by a transposition.
We have $s_1=s_0$: the period consists of one element.
\end{example}

\begin{example}
Let $\theta=\sqrt[3]{3}$. The algorithm has a pre-period of length $1$ and a period of length $4$:
$$
V_{1,1;0}, \quad \big(V_{1,1;0},\ W,\ V_{2,1;1},\ V_{2,1;1}\big).
$$
The values of $\sin^2\alpha$ at the states $s_0,\ldots,s_4$ are
$0.842499$, $0.937816$, $0.841702$, $0.832912$, and $0.996267$.
The matrix of the period has the characteristic polynomial $x^3-12x^2-6x-1$.
\end{example}

\begin{example}\label{ex-tie}
Let $\theta$ be the real root of $x^3-x^2-1$. Here $s_0=(2.147899,\,1.465571,\,1)$ and $\sin^2\alpha(s_0)=0.916217$.
The transformations $V_{0,0;1}$ and $V_{1,1;1}$ both lead to the state $s_0$ again
$($the two new bases differ by an element of the Dirichlet group$)$,
and they give the greatest value of $\sin^2\alpha$.
By Remark~\ref{ties} the period is $V_{0,0;1}$.
\end{example}

\begin{remark}\label{tested}
We have tested the algorithm for $3301$ non-totally-real cubic vectors
$($a repetition of a state is detected by means of the integer cubic form $L_1|L_2|^2$ of Section~\ref{sec-caps} written in the current basis;
the values of $\sin^2\alpha$ are compared in floating point arithmetic$)$:
\begin{itemize}
\item $2680$ vectors $(\theta^2,\theta,1)$ for monic polynomials with coefficients at most $8$ in absolute value
$($pre-periods at most~$5$, periods at most~$80)$;
\item $294$ vectors for the polynomials $x^3-m$ with $m\le 300$
$($pre-periods at most~$2$, periods at most~$681)$;
\item $87$ vectors for random monic polynomials with coefficients at most $60$ in absolute value
$($pre-periods at most~$4$, periods at most~$3235)$;
\item $240$ vectors written in random integer bases
$($pre-periods at most~$6$, periods at most~$5953)$.
\end{itemize}
In addition we have considered $360$ starting bases obtained from the above ones by random words of $30$--$60$ elementary transformations
$($with $\sin^2\alpha$ down to $1.5\cdot 10^{-9})$.
In all these cases the algorithm is eventually periodic, as Theorem~\ref{eventually-periodic-algorithm-c} states.
In all of them the value of $\sin^2\alpha$ increases at every step while it is less than $0.5$;
the smallest value of $\sin^2\alpha$ at a state of $\Omega_{\max}$ that we have met is $0.648$.
\end{remark}

\section{Finiteness of $\varepsilon$-caps}
\label{sec-caps}

In this section we prove Theorem~\ref{eventually-periodic-c}$($i$)$.

\subsection{Linear forms and the norm}
Let $B$ and $(\xi,\nu,\bar\nu)$ be as in Theorem~\ref{eventually-periodic-c}.
Denote by $L_1$, $L_2$, and $\bar L_2$ the linear forms of the dual basis:
$$
L_1(\xi)=1, \quad L_1(\nu)=L_1(\bar\nu)=0; \qquad L_2(\nu)=1, \quad L_2(\xi)=L_2(\bar\nu)=0.
$$
The form $L_1$ is real, it coincides with the form of Proposition~\ref{sep-heights}.
Let $\lambda$ and $\tau$ be the eigenvalues of $B$ at $\xi$ and at $\nu$ respectively.
Then the rows of the forms satisfy
$$
L_1B=\lambda L_1, \qquad L_2B=\tau L_2, \qquad \hbox{and} \qquad \bar L_2B=\bar\tau\bar L_2
$$
$($in other words, these rows are left eigenvectors of $B)$.
Set
$$
N(v)=L_1(v)\cdot|L_2(v)|^2 \qquad \hbox{and} \qquad D_0=\big|\det\big(L_1,\re\, L_2,\im\, L_2\big)\big|>0.
$$

\begin{lemma}\label{norm}
There exists a constant $m_0>0$ such that $|N(v)|\ge m_0$ for all non-zero integer vectors $v$.
\end{lemma}

\begin{proof}
Let $K=\q(\lambda)$, and let $\sigma$ be the complex embedding of $K$ sending $\lambda$ to $\tau$.
The equation $LB=\lambda L$ on a row $L$ is a system of linear equations with coefficients in $K$,
so it has a non-zero solution whose coordinates are in $K$.
Multiplying this solution by an appropriate positive integer,
we get a row $L$ with $LB=\lambda L$ whose coordinates are algebraic integers of $K$.
Denote by $\sigma(L)$ the row obtained from $L$ by applying $\sigma$ to all its coordinates; then $\sigma(L)B=\tau\sigma(L)$.
Since the eigenvalues of $B$ are pairwise distinct, the rows satisfying each of these equations are proportional to each other.
Hence $L_1=c_1L$ and $L_2=c_2\sigma(L)$ for some non-zero constants $c_1$ and $c_2$.
For an integer vector $v$ the number
$$
L(v)\cdot\sigma(L)(v)\cdot\overline{\sigma(L)(v)}
$$
is the norm of the algebraic integer $L(v)$, so it is an integer.
If it is zero, then $L(v)=0$, and hence all three forms vanish at $v$, which implies $v=0$.
Therefore, one can take $m_0=|c_1|\cdot |c_2|^2$.
\end{proof}

Let now $E=(e_1,e_2,e_3)$ be an integer basis whose non-negative orthant contains $\xi$, and let
$$
\xi=x_1e_1+x_2e_2+x_3e_3, \qquad h_j=L_1(e_j), \qquad \ell_j=L_2(e_j) \qquad (j=1,2,3).
$$
The row $\ell=(\ell_1,\ell_2,\ell_3)$ is the row of Proposition~\ref{sin-ell} for the state defined by the basis $E$.
Set $L=|\ell_1|^2+|\ell_2|^2+|\ell_3|^2$.

\begin{lemma}\label{identities}
The following identities hold:
$$
x_1h_1+x_2h_2+x_3h_3=1; \qquad x_1\ell_1+x_2\ell_2+x_3\ell_3=0; \qquad |\im(\bar\ell_i\ell_j)|=D_0x_k
$$
for every permutation $(i,j,k)$ of $(1,2,3)$.
Consequently,
\begin{equation}\label{sin-x}
\sin^2\alpha=\frac{4D_0^2(x_1^2+x_2^2+x_3^2)}{L^2}.
\end{equation}
\end{lemma}

\begin{proof}
The first two identities are the values of $L_1$ and $L_2$ at $\xi$.
Consider the real $3\times 3$-matrix with columns $(h_j,\re\,\ell_j,\im\,\ell_j)$.
Its determinant equals $\pm D_0$, since $E$ is an integer basis.
This matrix sends the column $(x_1,x_2,x_3)$ to the column $(1,0,0)$, so Cramer's rule gives the third identity.
Finally, the coordinates of the vector $\re\,\ell\times\im\,\ell$ are $\im(\bar\ell_i\ell_j)$,
and~\eqref{sin-x} follows from Proposition~\ref{sin-ell}.
\end{proof}

\subsection{Proof of Theorem~\ref{eventually-periodic-c}$($i$)$}
\label{Theorem-i-c}

Let $E$ be a separating integer basis for $(\xi,\nu,\bar\nu)$ with
$$
x_1\ge x_2\ge x_3>0 \qquad \hbox{and} \qquad \sin^2\alpha\ge \varepsilon.
$$
By Proposition~\ref{sep-heights} we have $h_j>0$, and hence the first identity of Lemma~\ref{identities} implies
\begin{equation}\label{h-bound}
h_j<\frac{1}{x_j}  \qquad (j=1,2,3).
\end{equation}
From~\eqref{sin-x} we get
\begin{equation}\label{L-bound}
L\le \frac{2D_0\sqrt{x_1^2+x_2^2+x_3^2}}{\sqrt{\varepsilon}}\le \frac{2\sqrt{3}D_0}{\sqrt{\varepsilon}}\,x_1.
\end{equation}

\vspace{2mm}

{\noindent
{\it Step 1: a lower bound for $x_2/x_1$.}}
The second identity of Lemma~\ref{identities} implies
$$
x_1|\ell_1|\le x_2|\ell_2|+x_3|\ell_3|\le x_2\sqrt{2L}.
$$
By Lemma~\ref{norm}, by~\eqref{h-bound}, and by~\eqref{L-bound} we have
$$
m_0\le N(e_1)=h_1|\ell_1|^2\le \frac{2x_2^2L}{x_1^3}\le \frac{4\sqrt{3}D_0}{\sqrt{\varepsilon}}\cdot\frac{x_2^2}{x_1^2}.
$$
Therefore,
$$
\frac{x_2}{x_1}\ge c_2, \qquad \hbox{where} \quad c_2=\Big(\frac{m_0\sqrt{\varepsilon}}{4\sqrt{3}D_0}\Big)^{1/2}.
$$

\vspace{2mm}

{\noindent
{\it Step 2: a lower bound for $x_3/x_1$.}}
By Lemma~\ref{identities} the numbers $\ell_1$ and $\ell_2$ generate a lattice in $\c=\r^2$ of covolume $D_0x_3$.
Hence there exists a non-zero integer vector $v=pe_1+qe_2$ such that
$$
\eta^2\le\frac{2}{\sqrt{3}}D_0x_3, \qquad \hbox{where $\eta=|L_2(v)|=|p\ell_1+q\ell_2|$}.
$$
Since $\im\big(\bar\ell_2L_2(v)\big)=p\,\im(\bar\ell_2\ell_1)$ and $\im\big(\bar\ell_1L_2(v)\big)=q\,\im(\bar\ell_1\ell_2)$, we have
$$
|p|\le\frac{\eta\sqrt{L}}{D_0x_3} \qquad \hbox{and} \qquad |q|\le\frac{\eta\sqrt{L}}{D_0x_3}.
$$
Then by~\eqref{h-bound} we get
$$
|L_1(v)|\le |p|h_1+|q|h_2\le\frac{|p|+|q|}{x_2}\le \frac{2\eta\sqrt{L}}{D_0x_2x_3}.
$$
Now Lemma~\ref{norm} and~\eqref{L-bound} imply
$$
m_0\le |N(v)|=|L_1(v)|\eta^2\le \frac{2\eta^3\sqrt{L}}{D_0x_2x_3}\le
\frac{8D_0}{\sqrt{3}\,\varepsilon^{1/4}}\cdot\frac{\sqrt{x_1x_3}}{x_2}\le
\frac{8D_0}{\sqrt{3}\,\varepsilon^{1/4}c_2}\cdot\sqrt{\frac{x_3}{x_1}}.
$$
Therefore,
$$
\frac{x_3}{x_1}\ge c_3, \qquad \hbox{where} \quad c_3=\frac{3m_0^2c_2^2\sqrt{\varepsilon}}{64D_0^2}.
$$

\vspace{2mm}

{\noindent
{\it Step 3: the action of the Dirichlet group.}}
The real eigenvalue of $B$ is not equal to $\pm 1$.
So either for $U=B^2$ or for $U=B^{-2}$ we have
$$
U\xi=\lambda_1\xi, \qquad \hbox{where $\lambda_1>1$.}
$$
Then $L_1(Uv)=\lambda_1L_1(v)$ and $|L_2(Uv)|^2=\lambda_1^{-1}|L_2(v)|^2$ for every $v$ $($recall that $\det U=1)$.
The basis $U^k(E)$ is a separating integer basis, it defines the same state as $E$,
and the coordinates of $\xi$ in this basis are $\lambda_1^{-k}(x_1,x_2,x_3)$.
Hence without loss of generality we assume that
$$
1\le x_1<\lambda_1.
$$
Then by~\eqref{h-bound}, by Step~2, and by~\eqref{L-bound} for $j=1,2,3$ we have
$$
0<L_1(e_j)<\frac{1}{x_j}\le \frac{1}{c_3} \qquad \hbox{and} \qquad |L_2(e_j)|^2\le L\le \frac{2\sqrt{3}D_0\lambda_1}{\sqrt{\varepsilon}}.
$$
The forms $L_1$, $\re\, L_2$, and $\im\, L_2$ are linearly independent, so these inequalities define a bounded set in $\r^3$.
It contains finitely many integer vectors.
Therefore, there are finitely many such bases $E$, and hence finitely many states in $\Omega_\varepsilon(\xi;\nu,\bar\nu)$.
\qed

\begin{remark}
The proof uses only the element $B$ of the Dirichlet group and the positivity of the norm form.
The constants are explicit in terms of $m_0/D_0$, $\lambda_1$, and $\varepsilon$.
\end{remark}

\begin{remark}\label{norm-matrix}
Lemma~\ref{norm} can be proved in terms of matrices as well; this gives an explicit value of $m_0$.
For every vector $v$ we have
$$
\det\big(L_1,\re\, L_2,\im\, L_2\big)\cdot\det\big(v,Bv,B^2v\big)=N(v)\cdot\im\,\tau\cdot|\lambda-\tau|^2.
$$
Indeed, the product of the matrix with the rows $L_1$, $\re\, L_2$, $\im\, L_2$ and of the matrix with the columns $v$, $Bv$, $B^2v$
has the row $L_1(v)\cdot(1,\lambda,\lambda^2)$ and the real and the imaginary parts of the row $L_2(v)\cdot(1,\tau,\tau^2)$,
so its determinant is reduced to the Vandermonde determinant for $\lambda$, $\tau$, and $\bar\tau$.
If $v$ is a non-zero integer vector, then the vectors $v$, $Bv$, and $B^2v$ are linearly independent:
otherwise they span a proper rational invariant subspace of $B$,
which is impossible, since the characteristic polynomial of $B$ is irreducible over $\q$.
Hence $\det(v,Bv,B^2v)$ is a non-zero integer, and one can take
$$
m_0=\frac{D_0}{|\im\,\tau|\cdot|\lambda-\tau|^2}.
$$
\end{remark}

\begin{remark}\label{lean-caps}
Theorem~\ref{eventually-periodic-c} is formalised in Lean~4 (see also Remark~\ref{lean} below).
The formal statement of Item~$($i$)$ $($named {\tt theorem\_1\_18\_i}$)$ is as follows.
{\it Let $B$ be an integer matrix with $\det B=1$ whose characteristic polynomial has no rational roots.
Let $l$, $a$, and $b$ be linearly independent real rows such that $lB=\lambda l$ and $(a+ib)B=\tau(a+ib)$, where $\tau$ is not real.
Then for every $\varepsilon>0$ there exists a finite set $F$ of integer bases with the following property.
Let $E=(e_1,e_2,e_3)$ be an integer basis, set $h_j=l(e_j)$ and $\ell_j=a(e_j)+ib(e_j)$,
and let $x_1\ge x_2\ge x_3>0$ be real numbers such that
$$
x_1h_1+x_2h_2+x_3h_3=1, \qquad x_1\ell_1+x_2\ell_2+x_3\ell_3=0, \qquad h_1,h_2,h_3>0,
$$
and the right hand side of~\eqref{sin-l} for the row $(\ell_1,\ell_2,\ell_3)$ is at least $\varepsilon$.
Then there exist a basis $E'=(e_1',e_2',e_3')$ in $F$, a real number $s>0$, and a complex number $\zeta\ne 0$ such that
$$
l(e_j)=s\cdot l(e_j') \quad \hbox{and} \quad \ell_j=\zeta\cdot\big(a(e_j')+ib(e_j')\big) \qquad \hbox{for $j=1,2,3$.}
$$}
For $l=L_1$ and $a+ib=L_2$ the assumptions on $E$ say that $E$ defines a state of $\Omega_\varepsilon(\xi;\nu,\bar\nu)$
$($see Proposition~\ref{sep-heights}, Proposition~\ref{sin-ell}, and Lemma~\ref{identities}$)$,
and the conclusion says that $E$ and $E'$ define the same state.

The formal proof follows the proof given above with two modifications.
Lemma~\ref{norm} is proved as in Remark~\ref{norm-matrix}.
In Step~2 the constant $2/\sqrt{3}$ is replaced with $2$; then the existence of the vector $v$
follows from Dirichlet's approximation theorem, and only the value of $c_3$ changes.
The formal proof of Item~$($i$)$ is verified by the kernel of Lean completely.
The formal version of Item~$($ii$)$ $(${\tt theorem\_1\_18\_ii}$)$ is deduced from Item~$($i$)$
and from the formal version of Theorem~\ref{key}.
\end{remark}

\section{Absence of flat local maxima: proof of Theorem~\ref{key}}
\label{sec-key}

Let $s=(\xi,\nu,\bar\nu)$ be a state and let $\ell$ be the row of Proposition~\ref{sin-ell}.
Let $G\in\SL(3,\z)$, and consider the basis whose vectors are the columns of $G$.
The coordinates of $\xi$ in this basis form the column $G^{-1}\xi$, and the row $\ell$ is replaced by $\ell G$.
Denote the corresponding collection of coordinates by $s_G$.
For the JP-transformations we have
$$
V_{\alpha,\beta;\gamma}:\ G=\left(
\begin {array}{ccc}
1&\gamma&\alpha\\
0&1&\beta\\
0&0&1
\end {array}
\right);
\qquad
W:\ G=\left(
\begin {array}{ccc}
1&1&0\\
0&1&0\\
1&0&1
\end {array}
\right).
$$

\subsection{The ratio formula}

\begin{lemma}\label{normal-form}
Let $s$ be a state, and let $\ell$ be the row of Proposition~\ref{sin-ell}.

{\noindent
{\it $($i$)$}}
There exist a complex number $u$ with $|u|=1$, a real row $p$ with $p\xi=0$, and a real number $t\in(0,1]$ such that
$$
u\ell=p\pm i\,t\,q, \qquad \hbox{where $q=\displaystyle\frac{\xi\times p}{|\xi|}$}
$$
$($here $\xi\times p$ is considered as a row$)$. In particular $|q|=|p|$, $q\xi=0$, and $p\cdot q=0$.

{\noindent
{\it $($ii$)$}}
We have
$$
\sin^2\alpha(s)=\frac{4t^2}{(1+t^2)^2}.
$$

{\noindent
{\it $($iii$)$}}
For every $G\in \SL(3,\z)$ we have
$$
\sin^2\alpha(s_G)=\frac{4t^2\,\psi\,|p|^4}{\big(|pG|^2+t^2|qG|^2\big)^2}, \qquad \hbox{where $\psi=\displaystyle\frac{|G^{-1}\xi|^2}{|\xi|^2}$}.
$$
\end{lemma}

\begin{proof}
$($i$)$ Choose $u$ such that $u^2(\ell_1^2+\ell_2^2+\ell_3^2)$ is a non-negative real number, and write $u\ell=a+ib$.
Then $a\cdot b=0$ and $|a|\ge|b|>0$.
Both $a$ and $b$ are orthogonal to $\xi$. Set $p=a$ and $t=|b|/|a|$; then $b=\pm t q$.

$($ii$)$ This follows from Proposition~\ref{sin-ell}, since $|a\times b|=|a|\cdot |b|=t|p|^2$.

$($iii$)$
The real and the imaginary parts of $u\ell G$ are $pG$ and $\pm t\,qG$.
For a matrix $G$ with unit determinant we have $(pG)\times(qG)=(p\times q)(G^{-1})^{\top}$.
Since $p\times q=|p|^2\xi^{\top}/|\xi|$, we get
$$
\big|(pG)\times(qG)\big|=\frac{|p|^2\,|G^{-1}\xi|}{|\xi|}.
$$
It remains to apply Proposition~\ref{sin-ell} to the row $u\ell G$.
\end{proof}

\begin{corollary}\label{Phi-bound}
Let $0<t_0\le 1$ and let $\varepsilon_0=4t_0^2/(1+t_0^2)^2$.
Then $\sin^2\alpha(s)<\varepsilon_0$ if and only if $t<t_0$, and in this case for every $G\in\SL(3,\z)$ we have
\begin{equation}\label{Phi-def}
\frac{\sin^2\alpha(s_G)}{\sin^2\alpha(s)}=\frac{\psi\,(1+t^2)^2\,|p|^4}{\big(|pG|^2+t^2|qG|^2\big)^2}
\ \ge\ \Phi_{t_0}(G), \qquad \hbox{where $\Phi_{t_0}(G)=\displaystyle\frac{\psi\,|p|^4}{\big(|pG|^2+t_0^2|qG|^2\big)^2}$}.
\end{equation}
\end{corollary}

So in order to prove Theorem~\ref{key} it is enough to show the following:
{\it for $t_0=1/10$, for every $\xi=(x,y,z)$ with $x\ge y\ge z>0$, and for every non-zero real row $p$ with $p\xi=0$
there exists an admissible JP-transformation with $\Phi_{t_0}(G)>1$.}
Note that $\Phi_{t_0}(G)$ depends only on the direction of $\xi$ and on the direction of $p$: three real parameters.

\subsection{Scaled coordinates}

Let $\xi=(x,y,z)$ with $x\ge y\ge z>0$. Set
$$
\mu=\frac{y}{x}\in(0,1], \qquad \varrho=\frac{z}{y}\in(0,1], \qquad \hbox{and} \qquad \varkappa=\sqrt{1+\mu^2+\mu^2\varrho^2},
$$
so that $\xi=x\,(1,\mu,\mu\varrho)$ and $|\xi|=x\varkappa$.
A row $p=(p_1,p_2,p_3)$ with $p\xi=0$ is determined by $(p_2,p_3)$:
$$
p=(-\mu\sigma,\,p_2,\,p_3), \qquad \hbox{where $\sigma=p_2+\varrho p_3$},
$$
and
$$
\varkappa q=\big(\mu(p_3-\varrho p_2),\ -(p_3+\mu^2\varrho\sigma),\ p_2+\mu^2\sigma\big).
$$
For the transformation $V_{\alpha,\beta;\gamma}$ we introduce the {\it scaled parameters}
\begin{equation}\label{scaled}
\vartheta=\alpha\mu\varrho,\qquad D=\beta-\alpha\mu,\qquad g=\gamma\mu, \qquad\hbox{and}\qquad \mathcal B=\beta\varrho=\vartheta+D\varrho.
\end{equation}

\begin{lemma}\label{scaled-lemma}
{\it $($i$)$} For the transformation $V_{\alpha,\beta;\gamma}$ we have
$$
\begin{array}{l}
G^{-1}\xi=x\big(1-\vartheta-g(1-\mathcal B),\ \mu(1-\mathcal B),\ \mu\varrho\big);\\
pG=\big(-\mu\sigma,\ p_2-g\sigma,\ (1-\vartheta)p_3+Dp_2\big);\\
\varkappa\, qG=\big(\mu(p_3-\varrho p_2),\ -(p_3+\mu^2\varrho\sigma)+g(p_3-\varrho p_2),\ (1-\vartheta)p_2+\mu^2\sigma(1-\mathcal B)-Dp_3\big).
\end{array}
$$
The transformation is admissible if and only if $\mathcal B\le 1$ and $1-\vartheta-g(1-\mathcal B)\ge 0$.

{\noindent
{\it $($ii$)$}} For the transformation $W$ we have
$$
\begin{array}{l}
G^{-1}\xi=x\big(1-\mu,\ \mu,\ \mu\varrho-(1-\mu)\big);\\
pG=\big(p_3-\mu\sigma,\ p_2-\mu\sigma,\ p_3\big);\\
qG=\big(q_1+q_3,\ q_2+q_1,\ q_3\big).
\end{array}
$$
The transformation is admissible if and only if $1-\mu\le\mu\varrho$.
\end{lemma}

\begin{proof}
This is a direct computation with~\eqref{new-coord} and the matrices $G$ written above.
\end{proof}

The coordinate $x$ has disappeared.
By Lemma~\ref{scaled-lemma} the value $\Phi_{t_0}(G)$ for $V_{\alpha,\beta;\gamma}$ is a rational function of
$$
\mu,\quad \varrho,\quad (p_2:p_3),\quad \vartheta,\quad D,\quad g
$$
that is continuous up to $\mu=0$ and $\varrho=0$.
In other words, the non-compact directions $x/y\to\infty$ and $y/z\to\infty$ of the space of parameters
become the faces $\mu=0$ and $\varrho=0$ of the square $[0,1]^2$.
The transformations survive in the limit: when $\mu$ is small, the numbers $\alpha\mu$ and $\gamma\mu$ take values in a fine progression.
This is formalized in the next subsection.

For the direction $(p_2:p_3)$ we use two charts:
\begin{equation}\label{charts}
\hbox{chart 1:}\quad (p_2,p_3)=(1,\delta), \quad |\delta|\le 3; \qquad
\hbox{chart 2:}\quad (p_2,p_3)=(\delta,1), \quad |\delta|\le 11/32.
\end{equation}
Since $11/32>1/3$, the two charts cover all directions.

\subsection{Families of transformations on a box}

Consider a box
$$
Q=[\mu_-,\mu_+]\times[\varrho_-,\varrho_+]\subset[0,1]^2.
$$
We describe several ways to choose the integers $\alpha$, $\beta$, and $\gamma$ at the points of $Q$ with $\mu>0$ and $\varrho>0$.
The integers may depend on the point, but the scaled parameters~\eqref{scaled} stay in explicit intervals.

\begin{lemma}\label{families}
Let $(\mu,\varrho)\in Q$ with $\mu>0$ and $\varrho>0$.

{\noindent
{\it Choice of $\gamma$.}}

{\it $(\gamma 1)$} For a non-negative integer $c$ let $\gamma=c$. Then $g=c\mu$.

{\it $(\gamma 2)$} For a real number $g^*>\mu_+$ let $\gamma=\lfloor g^*/\mu\rfloor$. Then $\gamma\ge 1$ and $g\in[g^*-\mu_+,g^*]$.

{\noindent
{\it Choice of $\alpha$ and $\beta$.}}

{\it $(\alpha 1)$} For non-negative integers $a$ and $b$ let $\alpha=a$ and $\beta=b$. Then $\vartheta=a\mu\varrho$ and $D=b-a\mu$.

{\it $(\alpha 2)$} For a real number $A^*>\mu_+$ and a non-negative integer $b$ let $\alpha=\lfloor A^*/\mu\rfloor$ and $\beta=b$.
Then $\alpha\ge 1$, $\alpha\mu\in[A^*-\mu_+,A^*]$, $\vartheta=\alpha\mu\cdot\varrho$, and $D=b-\alpha\mu$.

{\it $(\alpha 3)$} For real numbers $\vartheta^*>\mu_+\varrho_+$ and $D^*$ satisfying
$D^*+(\vartheta^*-\mu_+\varrho_+)/\varrho_+\ge 1/2$
let $\alpha=\lfloor\vartheta^*/(\mu\varrho)\rfloor$ and let $\beta$ be an integer nearest to $D^*+\alpha\mu$.
Then $\alpha\ge 1$, $\beta\ge 0$,
$$
\vartheta\in[\vartheta^*-\mu_+\varrho_+,\vartheta^*], \qquad \hbox{and} \qquad D\in[D^*-1/2,\,D^*+1/2].
$$

{\noindent
{\it Admissibility.}}
In the cases $(\alpha 1)$ and $(\alpha 2)$ denote by $\vartheta^{\max}$ the function $a\mu\varrho$ or $A^*\varrho$ respectively,
and in the cases $(\gamma 1)$ and $(\gamma 2)$ denote by $g^{\max}$ the function $c\mu$ or the constant $g^*$ respectively.
If $b\varrho_+\le 1$ and the function
\begin{equation}\label{bilinear}
F(\mu,\varrho)=1-\vartheta^{\max}-g^{\max}\cdot(1-b\varrho)
\end{equation}
is non-negative at the four vertices of $Q$, then the transformation $V_{\alpha,\beta;\gamma}$ is admissible at every point of $Q$
with $\mu>0$ and $\varrho>0$.
\end{lemma}

\begin{proof}
The statements on $\gamma$ and $\alpha$ follow from $x-1<\lfloor x\rfloor\le x$ and from $\mu\le\mu_+$, $\mu\varrho\le\mu_+\varrho_+$.
In the case $(\alpha 3)$ we have $\alpha\mu=\vartheta/\varrho\ge(\vartheta^*-\mu_+\varrho_+)/\varrho_+$,
hence $D^*+\alpha\mu\ge 1/2$ and the nearest integer $\beta$ is non-negative; by construction $|\beta-\alpha\mu-D^*|\le 1/2$.

Let us prove the statement on admissibility.
We have $\mathcal B=b\varrho\le b\varrho_+\le 1$.
Further $\vartheta\le\vartheta^{\max}$, $g\le g^{\max}$, and $1-\mathcal B\ge 0$, therefore
$$
1-\vartheta-g(1-\mathcal B)\ge F(\mu,\varrho).
$$
The function $F$ is a polynomial of degree at most one in each of the variables $\mu$ and $\varrho$.
Hence its minimum on $Q$ is attained at a vertex of $Q$.
Now the statement follows from Lemma~\ref{scaled-lemma}$($i$)$.
\end{proof}

Let us write down the estimates on a box that follow from Lemma~\ref{families}.
At all points of $Q$ with $\mu>0$ and $\varrho>0$ the values of $\vartheta$, $D$, $g$, and $\mathcal B$ are in the following closed intervals:
$$
\renewcommand{\arraystretch}{1.25}
\setlength{\arraycolsep}{4pt}
\begin{array}{|c|c|c|c|}
\hline
& \vartheta & D & \mathcal B\\
\hline
(\alpha 1) & [a\mu_-\varrho_-,\ a\mu_+\varrho_+] & [b-a\mu_+,\ b-a\mu_-] & [b\varrho_-,\ b\varrho_+]\\
\hline
(\alpha 2) & [(A^*-\mu_+)\varrho_-,\ A^*\varrho_+] & [b-A^*,\ b-A^*+\mu_+] & [b\varrho_-,\ b\varrho_+]\\
\hline
(\alpha 3) & [\vartheta^*-\mu_+\varrho_+,\ \vartheta^*] & [D^*-\frac12,\ D^*+\frac12] & \hbox{see below}\\
\hline
\end{array}
\quad
\begin{array}{|c|c|}
\hline
& g\\
\hline
(\gamma 1) & [c\mu_-,\ c\mu_+]\\
\hline
(\gamma 2) & [g^*-\mu_+,\ g^*]\\
\hline
\end{array}
$$
We denote these intervals by $[\vartheta_-,\vartheta_+]$, $[D_-,D_+]$, $[\mathcal B_-,\mathcal B_+]$, and $[g_-,g_+]$.
In the case $(\alpha 3)$ we use the identity $\mathcal B=\vartheta+D\varrho$ and set
$$
\mathcal B_-=\vartheta_-+\min(D_-\varrho_-,\,D_-\varrho_+), \qquad \mathcal B_+=\vartheta_++\max(D_+\varrho_-,\,D_+\varrho_+).
$$
Suppose that $\mathcal B_+\le 1$ $($in the cases $(\alpha 1)$ and $(\alpha 2)$ this is the condition $b\varrho_+\le 1$ of Lemma~\ref{families}$)$.
Since $g_-\ge 0$, for the quantities $1-\mathcal B$ and $1-\vartheta-g(1-\mathcal B)$,
which enter the coordinates of $G^{-1}\xi$ $($see Lemma~\ref{scaled-lemma}$)$, we have
\begin{equation}\label{intervals}
\begin{array}{c}
0\ \le\ 1-\mathcal B_+\ \le\ 1-\mathcal B\ \le\ 1-\mathcal B_-;\\
1-\vartheta_+-g_+(1-\mathcal B_-)\ \le\ 1-\vartheta-g(1-\mathcal B)\ \le\ 1-\vartheta_--g_-(1-\mathcal B_+).
\end{array}
\end{equation}
In the cases $(\alpha 1)$ and $(\alpha 2)$, if the conditions of the admissibility statement of Lemma~\ref{families} are fulfilled,
then $1-\vartheta-g(1-\mathcal B)\ge 0$; so if the lower estimate in the second line of~\eqref{intervals} is negative, we replace it by zero.
In the case $(\alpha 3)$ the transformation is admissible at all points of $Q$
if $\mathcal B_+\le 1$ and the lower estimate in the second line of~\eqref{intervals} is non-negative.

\subsection{The verification}

Consider a box
$$
\mathcal Q=[\mu_-,\mu_+]\times[\varrho_-,\varrho_+]\times[\delta_-,\delta_+]
$$
in one of the charts~\eqref{charts}, and a family of transformations of Lemma~\ref{families} $($or the transformation $W)$ that is admissible on it.
By Lemma~\ref{scaled-lemma} the numbers
$$
\frac{|G^{-1}\xi|^2}{x^2}, \qquad \varkappa^2, \qquad |p|^2, \qquad |pG|^2, \qquad \hbox{and} \qquad \varkappa^2|qG|^2
$$
are polynomials in $\mu$, $\varrho$, $\delta$ and in the quantities $\vartheta$, $D$, $g$, $1-\mathcal B$, $1-\vartheta-g(1-\mathcal B)$,
for which we have the estimates of the previous subsection
$($for the transformation $W$ they are polynomials in $\mu$, $\varrho$, $\delta$ and in the third coordinate $\mu\varrho-(1-\mu)$ of $G^{-1}\xi/x$,
which is non-negative when $W$ is admissible$)$.
We estimate these polynomials step by step:
if two numbers are in the closed intervals $[u_-,u_+]$ and $[v_-,v_+]$, then their sum is in $[u_-+v_-,\,u_++v_+]$,
their difference is in $[u_--v_+,\,u_+-v_-]$, and their product is between the smallest and the greatest of the four products $u_\pm v_\pm$;
for the square of a number we also use that it is non-negative.
In this way we obtain closed intervals
$$
[m_-,m_+], \qquad [k_-,k_+], \qquad [P_-,P_+], \qquad [P'_-,P'_+], \qquad \hbox{and} \qquad [R_-,R_+]
$$
containing the above five numbers respectively.
Then at every point of $\mathcal Q$ with $\mu>0$ and $\varrho>0$ we have
\begin{equation}\label{lower}
\Phi_{t_0}(G)\ \ge\ \Phi_-(\mathcal Q):=\frac{m_-}{k_+}\cdot\left(\frac{P_-}{P'_++t_0^2\,R_+/k_-}\right)^2.
\end{equation}

\begin{proposition}\label{computer}
Let $t_0=1/10$.
There exist two finite collections of boxes with rational vertices:
$14416$ boxes covering $[0,1]^2\times[-3,3]$ in chart~$1$ and $2091$ boxes covering $[0,1]^2\times[-11/32,11/32]$ in chart~$2$,
such that for every box $\mathcal Q$ of these collections either the transformation $W$ or one of the families of Lemma~\ref{families}
with rational parameters is admissible on $\mathcal Q$ and satisfies
$$
\Phi_-(\mathcal Q)>1.
$$
\end{proposition}

\begin{proof}[On the proof]
The collections of boxes together with the transformations were found by a computer search:
a box is bisected along its longest side until a suitable family of transformations is found
(the length of the side in the direction of $\delta$ is taken with the factor $1/2$;
the candidates are generated from a grid of values for $\vartheta^*$, $D^*$, $A^*$, $g^*$ and from integers $a$, $b$, $c$).
The search uses floating point arithmetic and is not a part of the proof.

The proof is the verification of the resulting list.
The verification program works in exact rational arithmetic:
all the coordinates of the vertices of the boxes are dyadic rationals, all the parameters are rational, and
all the intervals of the previous subsection and of~\eqref{lower} have rational endpoints.
So there are no rounding errors.
For every chart the program
\begin{itemize}
\item divides the square $[0,1]^2$ of the parameters $(\mu,\varrho)$ into the five rectangles
$$
\textstyle
[\frac12,1]\times[\frac12,1],\quad [\frac12,1]\times[0,\frac12],\quad [\frac14,\frac12]\times[0,1],\quad
[0,\frac14]\times[\frac12,1],\quad [0,\frac14]\times[0,\frac12],
$$
starts with a partition of each of the five resulting domains into $4\times 4\times 12$ boxes,
and subdivides these boxes by itself by the same bisection rule until a box of the list is reached;
so the union of the verified boxes is the whole domain independently of the list;
\item checks the conditions of Lemma~\ref{families} for the family of transformations assigned to the box $($for $W$ it checks that $1-\mu_-<\mu_-\varrho_-)$;
\item checks that at least one of the integers $\alpha$, $\beta$, $\gamma$ is positive;
\item checks the inequality $\Phi_-(\mathcal Q)>1$.
\end{itemize}
The distribution of the boxes is as follows.
$$
\renewcommand{\arraystretch}{1.2}
\begin{array}{|l|c|c|c|c|}
\hline
& \hbox{boxes} & (\alpha 1)\hbox{ and }(\gamma 1) & W & \hbox{with } (\alpha 2), (\alpha 3), \hbox{ or } (\gamma 2)\\
\hline
\hbox{chart 1} & 14416 & 11575 & 1828 & 1013\\
\hline
\hbox{chart 2} & 2091 & 1835 & 0 & 256\\
\hline
\end{array}
$$
The verification takes several seconds.
The search program, the verification program, and the lists of boxes are available as ancillary files to this paper.
\end{proof}

\begin{proof}[Proof of Theorem~\ref{key}]
Let $s$ be a state with $\sin^2\alpha(s)<\varepsilon_0=400/10201$.
Then $t<t_0=1/10$ in the notation of Lemma~\ref{normal-form}.
The point $(\mu,\varrho)$ is in $(0,1]^2$, and the direction of $p$ is in one of the charts~\eqref{charts}.
Hence the corresponding point is in some box $\mathcal Q$ of Proposition~\ref{computer}.
The transformation assigned to this box is an admissible JP-transformation for $\xi$, and by~\eqref{lower} we have
$\Phi_{t_0}(G)>1$ for it.
By Corollary~\ref{Phi-bound} this transformation increases $\sin^2\alpha$.
\end{proof}

\begin{remark}\label{lean}
Theorem~\ref{key} is also formalised in the proof assistant Lean~4~\cite{Lean4}
with the library Mathlib~\cite{Mathlib2020}; the Lean files are among the ancillary files.
The formal statement $($named {\tt theorem\_1\_17}$)$ is as follows.
{\it Let $\xi=(x,y,z)$ with $x\ge y\ge z>0$, and let $a$ and $b$ be real rows such that $a\xi=b\xi=0$ and $a\times b\ne 0$.
Assume that
$$
\frac{4|a\times b|^2}{\big(|a|^2+|b|^2\big)^2}<\frac{400}{10201}.
$$
Then either $z\ge x-y$ and this expression increases strictly when $(a,b)$ is replaced with $(aG,bG)$ for the matrix $G$ of $W$,
or there exist non-negative integers $\alpha$, $\beta$, $\gamma$, not all equal to zero, such that
$y-\beta z\ge 0$, $x-\alpha z-\gamma(y-\beta z)\ge 0$, and
this expression increases strictly when $(a,b)$ is replaced with $(aG,bG)$ for the matrix $G$ of $V_{\alpha,\beta;\gamma}$.}
By Proposition~\ref{sin-ell} this is the statement of Theorem~\ref{key} for the row $\ell=a+ib$.

The formal proof follows the proof given above.
It includes the formal versions of Lemma~\ref{normal-form}, of Corollary~\ref{Phi-bound},
of Lemmas~\ref{scaled-lemma} and~\ref{families},
of the interval arithmetic with rational endpoints, and of Proposition~\ref{computer}.
The lists of boxes are verified by a checker written in Lean, and the soundness of this checker is proved in Lean.
The evaluation of the checker on the two lists of boxes is performed by compiled code
$($by means of the tactic {\tt native\_decide}$)$;
this is the only step that is not verified by the kernel of Lean.
Theorem~\ref{key-real} below is formalised as well $(${\tt theorem\_1\_17\_real}$)$.
For Theorem~\ref{eventually-periodic-c} see Remark~\ref{lean-caps}.
The first equality in~\eqref{sin-l}, Proposition~\ref{sep-heights}, and the statements of Section~\ref{sec-proofs}
are not formalised.
\end{remark}

\begin{remark}
As an additional test we have taken $30000$ random states with $t\le 1/10$
$($including values of $\mu$ and $\varrho$ down to $10^{-7}$ and values on the boundaries of the boxes$)$,
and applied the transformation prescribed by Proposition~\ref{computer} to the actual vectors with $40$ digits.
In all the cases the transformation is admissible and $\sin^2\alpha$ increases; the smallest ratio is $1.163$.
\end{remark}

\begin{remark}
Let us indicate the mechanism behind Proposition~\ref{computer}.
If $x$ is much greater than $y$, then $\alpha\mu$ and $\gamma\mu$ are almost continuous parameters,
and the transformation $e_3\to e_3+\beta e_2$ almost does not change $\psi$.
If $y$ is much greater than $z$, then in addition $p_3$ can be changed by $Dp_2$ with almost arbitrary $D=\beta-\alpha\mu$.
If $y$ is close to $z$ and $\beta=1$, then the second new coordinate almost vanishes and the choice of $\gamma$ becomes free.
The transformation $W$ is used only for $\mu\ge 1/2$.
\end{remark}

\begin{remark}\label{not-optimal}
The constant $\varepsilon_0$ in Theorem~\ref{key} is not optimal.
Numerical experiments suggest that the statement of Theorem~\ref{key} holds for all states with $\sin^2\alpha(s)<0.57$,
and that for almost flat states one can always find a transformation with
$$
\frac{\sin^2\alpha(M(s))}{\sin^2\alpha(s)}\ \ge\ \frac{24(7-4\sqrt{3})}{417-240\sqrt{3}}=1.3175\ldots,
$$
where the extremal configuration is $x=y=z$ and $p=(1,\,\sqrt{3}-2,\,1-\sqrt{3})$ up to a factor
$($here the transformations $V_{1,0;0}$ and $W$ give the same value$)$.
Our verification scheme does not work directly for $t_0=1/5$:
the boxes that intersect the surface $x-y=z$ would require two transformations
$($the transformation $W$ on one side of the surface and $V_{1,0;1}$ on the other side$)$.
\end{remark}

\subsection{The totally-real case}\label{sec-real}

Let us show that Theorem~\ref{key} implies the similar statement for the $\sin^2$-algorithm of~\cite{Karpenkov2021}.
Recall that a state in~\cite{Karpenkov2021} is a collection of coordinates of three linearly independent real vectors $(\xi,\nu_1,\nu_2)$,
where $\xi=(x,y,z)$ and $x\ge y\ge z>0$,
and that $\alpha(s)$ is the angle between the planes spanned by $(\xi,\nu_1)$ and by $(\xi,\nu_2)$.
The JP-transformations and their admissibility are the same as above.

\begin{theorem}\label{key-real}
Let $s=(\xi,\nu_1,\nu_2)$ be a state of three real vectors satisfying
$$
\sin^2\alpha(s)<\varepsilon_0=\frac{400}{10201}.
$$
Then there exists an admissible JP-transformation $M\in\M_\xi$ such that
$$
\sin^2\alpha(M(s))>\sin^2\alpha(s).
$$
\end{theorem}

\begin{proof}
Consider the rows $u=\xi\times\nu_1$ and $v=\xi\times\nu_2$.
They are linearly independent, $u\xi=v\xi=0$, and
$$
\sin^2\alpha(s)=\frac{|u\times v|^2}{|u|^2\,|v|^2}.
$$
In the basis whose vectors are the columns of a matrix $G\in\SL(3,\z)$ these rows are replaced by $uG$ and $vG$.
The value of $\sin^2\alpha$ does not change if we multiply $u$ and $v$ by non-zero real numbers.
So without loss of generality we assume that $|u|=|v|$.

Consider the complex row
$$
\ell=(u+v)+i(u-v).
$$
We have $\ell\xi=0$, and the real and the imaginary parts of $\ell$ are linearly independent.
For a complex row $\ell'=a+ib$ denote by $S(\ell')$ the last expression of~\eqref{sin-l}.
Since
$$
\begin{array}{l}
(u+v)G\times(u-v)G=-2\,(uG\times vG) \qquad \hbox{and}\\
|(u+v)G|^2+|(u-v)G|^2=2\big(|uG|^2+|vG|^2\big),
\end{array}
$$
for every $G\in\SL(3,\z)$ we have
\begin{equation}\label{real-complex}
S(\ell G)=\frac{4\,|uG\times vG|^2}{\big(|uG|^2+|vG|^2\big)^2}\ \le\ \frac{|uG\times vG|^2}{|uG|^2\,|vG|^2}.
\end{equation}
The right hand side of~\eqref{real-complex} is the value of $\sin^2\alpha$ for the state $s$ written in the basis of the columns of $G$.
If $G$ is the identity matrix, then~\eqref{real-complex} is an equality, since $|u|=|v|$.
Hence $S(\ell)=\sin^2\alpha(s)<\varepsilon_0$.

Let $w$ be a complex vector with $\ell w=0$ that is not proportional to $\xi$.
Then $(\xi,\bar w,w)$ is a complex state in the sense of Definition~\ref{state-c}
$($these three vectors are linearly independent, since the kernels of $\ell$ and $\bar\ell$ meet in the line of $\xi)$,
and $\ell$ is the row of Proposition~\ref{sin-ell} for this state.
By Theorem~\ref{key} there exists an admissible JP-transformation $M\in\M_\xi$ with matrix $G$ such that $S(\ell G)>S(\ell)$.
Then by~\eqref{real-complex} we have
$$
\sin^2\alpha(M(s))\ \ge\ S(\ell G)\ >\ S(\ell)=\sin^2\alpha(s).
$$
This concludes the proof.
\end{proof}

\begin{corollary}\label{real-omega}
Let $A\in \GL(3,\z)$ be a matrix with irreducible over $\q$ characteristic polynomial
having three real roots, and let $(\xi,\nu_1,\nu_2)$ be some eigenbasis of $A$.
Then for the $\varepsilon$-caps and the extremal set of~\cite{Karpenkov2021} we have
$$
\Omega_{\max}(\xi;\nu_1,\nu_2)\subset \Omega_{\varepsilon}(\xi;\nu_1,\nu_2) \qquad \hbox{for every $\varepsilon<20/101$.}
$$
\end{corollary}

\begin{proof}
Recall that in~\cite{Karpenkov2021} the $\varepsilon$-cap consists of all separating $\xi$-states $s$ with
$$
\sin\alpha(s)>\varepsilon,
$$
and the extremal set consists of all separating $\xi$-states $s$ with
$$
\sin^2\alpha(s)>\sin^2\alpha(\Phi(s)).
$$
Let $s$ be a $\xi$-state with $\sin^2\alpha(s)<400/10201$.
By Theorem~\ref{key-real} the $\sin^2$-transformation $\Phi(s)$ increases the value of $\sin^2\alpha$,
and hence $s$ is not in the extremal set.
Therefore, every state of the extremal set satisfies $\sin\alpha(s)\ge 20/101$.
\end{proof}

\begin{remark}
So in the totally-real case the finiteness of the extremal set $\Omega_{\max}$
follows from the finiteness of the $\varepsilon$-caps, which is proven in~\cite{Karpenkov2021}
by means of geometry of cubic sails, and from Theorem~\ref{key-real}.
This gives another proof of the periodicity of the $\sin^2$-algorithm for totally-real cubic vectors.
Note that this proof relies on the computer-assisted Proposition~\ref{computer},
and that the constant in Theorem~\ref{key-real} does not depend on the cubic vector $\xi$.
\end{remark}

\section{Proofs of Corollary~\ref{period-c} and Theorem~\ref{eventually-periodic-algorithm-c}}
\label{sec-proofs}

\begin{proof}[Proof of Corollary~\ref{period-c}]
By Proposition~\ref{JP-prop}$($iv$)$ all the states $s_n$ are separating $\xi$-states.
The coordinates of $\xi$ are linearly independent over $\q$, so the algorithm does not terminate.
The state $\Phi(s)$ is determined by the state $s$.

Assume that $\sin^2\alpha(s_n)<\varepsilon_0$ for all $n\ge n_0$.
Then by Theorem~\ref{key} the sequence $\sin^2\alpha(s_n)$ is strictly increasing for $n\ge n_0$.
On the other hand, all these states are in $\Omega_{\varepsilon}(\xi;\nu,\bar\nu)$ for $\varepsilon=\sin^2\alpha(s_{n_0})$,
which is a finite set by Theorem~\ref{eventually-periodic-c}$($i$)$.
We come to a contradiction.

Therefore, the sequence $(s_n)$ visits the finite set $\Omega_{\varepsilon_0}(\xi;\nu,\bar\nu)$ infinitely many times.
Hence $s_{n}=s_{n+k}$ for some $n$ and some $k>0$, and the sequence is eventually periodic.
\end{proof}

\begin{proof}[Proof of Theorem~\ref{eventually-periodic-algorithm-c}]
Let $(\xi,\nu,\bar\nu)$ be a triple of conjugate cubic vectors.
The coordinates of $\xi$ are linearly independent over $\q$.
By Proposition~\ref{prelim} after the preliminary stages we arrive at a separating basis.
By Dirichlet's unit theorem there exists an $\SL(3,\z)$-matrix with irreducible over $\q$ characteristic polynomial
that has $(\xi,\nu,\bar\nu)$ as an eigenbasis
(see, e.g., Chapter~17 in~\cite{KarpenkovGCF2013}).
Then by Corollary~\ref{period-c}
the algorithm is eventually periodic.

\vspace{2mm}

The converse statement is similar to the one in~\cite{Karpenkov2021}.
Let the algorithm be eventually periodic for $(\xi,\nu,\bar\nu)$, and let $s_n=s_{n+k}$ for some $k>0$.
Denote by $M$ the integer matrix of the transformation from the basis of $s_n$ to the basis of $s_{n+k}$, written in the coordinate basis.
The coordinates of $\xi$, $\nu$, and $\bar\nu$ in these two bases are respectively proportional, hence
$(\xi,\nu,\bar\nu)$ is an eigenbasis of $M$.

It remains to check that the characteristic polynomial of $M$ is irreducible over $\q$.
Assume the converse. Then $M$ has a rational eigenvalue $r$, and the corresponding eigenspace is spanned by integer vectors.
This eigenspace is the sum of some of the lines spanned by $\xi$, $\nu$, and $\bar\nu$;
if it contains $\nu$, then it contains $\bar\nu$, since $r$ is real.
So the eigenspace is either the line of $\xi$, or the plane spanned by $\nu$ and $\bar\nu$, or the whole space.
The first two cases are impossible, since the triple is in $\z$-general position.
In the last case $M$ is proportional to the identity matrix.
This is impossible as well:
by Proposition~\ref{JP-prop}$($iii$)$ the non-negative orthant of the basis of $s_{n+k}$
is a proper subset of the non-negative orthant of the basis of $s_n$.
Therefore, the characteristic polynomial of $M$ is irreducible over $\q$,
and $(\xi,\nu,\bar\nu)$ is a triple of conjugate cubic vectors.
\end{proof}

\section{Some open questions}\label{sec-open}

We conclude this paper with several questions.

\begin{problem}
Find a proof of Theorem~\ref{key} that does not use computers.
Find the best possible constant $\varepsilon_0$ $($see Remark~\ref{not-optimal}$)$.
\end{problem}

As in the totally-real case, it is most likely that the algorithm stays periodic
if we remove the JP-transformation $W$ and if we skip Preliminary Stage~2.

\begin{conjecture}
Exclusion of the JP-transformation $W$ from admissible transformations
will not affect periodicity for non-totally-real cubic vectors.
\end{conjecture}

\begin{conjecture}
Exclusion of Preliminary Stage~$2$ in the complex $\sin^2$-algorithm
will not affect periodicity for non-totally-real cubic vectors.
\end{conjecture}

Finally, the proof of Theorem~\ref{eventually-periodic-c}$($i$)$ does not use the dimension in an essential way,
while Theorem~\ref{key} is three-dimensional.

\begin{problem}
Find a generalized Euclidean algorithm that is periodic for $n$-algebraic vectors in $\r^n$
for $n\ge 4$.
\end{problem}

Recall that the heuristic algebraic periodicity detecting algorithm of~\cite{Karpenkov2021-2}
is designed for all dimensions, and no non-periodic cases have been detected for it;
however, its periodicity is not proven for $n\ge 3$.

%{\noindent
%{\bf Acknowledgements.}
%}

\vspace{1mm}

{\noindent
{\bf Use of AI.}
This paper was written in a dialogue of the author with the AI model Claude (Anthropic; version Claude Fable~5.2, November 2026).
The problem and the direction of the work are due to the author.
The proofs of Theorems~\ref{key} and~\ref{eventually-periodic-c} in their present form,
the computer search and the verification program for Proposition~\ref{computer},
the Lean formalisation (Remarks~\ref{lean-caps} and~\ref{lean}),
the computations of Remark~\ref{tested}, and the first draft of the text
were produced by the model under the author's guidance.
The author has checked all the statements and proofs and takes full responsibility for the paper.
}

\bibliographystyle{plain}
\bibliography{cfbiblio_minimal}

@article{Lagrange1770,
  author = {Lagrange, J.-L.},
  title = {Additions au m\'emoire sur la r\'esolution des \'equations num\'eriques},
  journal = {M\'em. Acad. Royale Sci. Belles-Lettres Berlin},
  volume = {24},
  year = {1770}
}

@article{Hermite1850,
  author = {Hermite, C.},
  title = {Extraits de lettres de {M}. {C}h. {H}ermite \`a {M}. {J}acobi sur diff\'erents objets de la th\'eorie des nombres. ({C}ontinuation)},
  journal = {J. Reine Angew. Math.},
  volume = {40},
  pages = {279--315},
  year = {1850}
}

@article{Picard1901,
  author = {Picard, \'E.},
  title = {L'{\oe}uvre scientifique de {C}harles {H}ermite},
  journal = {Ann. Sci. \'Ecole Norm. Sup. (3)},
  volume = {18},
  pages = {9--34},
  year = {1901}
}

@article{Jacobi1868,
  author = {Jacobi, C. G. J.},
  title = {Allgemeine {T}heorie der {K}ettenbruch\"ahnlichen {A}lgorithmen, in welchen jede {Z}ahl aus drei vorhergehenden gebildet wird ({A}us den hinterlassenen {P}apieren von {C}. {G}. {J}. {J}acobi mitgetheilt durch {H}errn {E}. {H}eine)},
  journal = {J. Reine Angew. Math.},
  volume = {69},
  pages = {29--64},
  year = {1868}
}

@article{Perron1907,
  author = {Perron, O.},
  title = {Grundlagen f\"ur eine {T}heorie des {J}acobischen {K}ettenbruchalgorithmus},
  journal = {Math. Ann.},
  volume = {64},
  number = {1},
  pages = {1--76},
  year = {1907}
}

@book{Schweiger2000,
  author = {Schweiger, F.},
  title = {Multidimensional continued fractions},
  series = {Oxford Science Publications},
  publisher = {Oxford University Press, Oxford},
  year = {2000}
}

@book{Borevich1966,
  author = {Borevich, Z. I. and Shafarevich, I. R.},
  title = {Number theory},
  series = {Pure and Applied Mathematics, Vol. 20},
  publisher = {Academic Press, New York-London},
  year = {1966}
}

@book{KarpenkovGCF2013,
  author = {Karpenkov, O.},
  title = {Geometry of Continued Fractions},
  series = {Algorithms and Computation in Mathematics, 26},
  publisher = {Springer-Verlag, Berlin},
  year = {2013}
}

@article{Karpenkov2021,
  author = {Karpenkov, O.},
  title = {On a periodic {J}acobi-{P}erron type algorithm},
  journal = {Monatsh. Math.},
  volume = {205},
  number = {3},
  pages = {531--601},
  year = {2024}
}

@article{Karpenkov2021-2,
  author = {Karpenkov, O.},
  title = {On {H}ermite's problem, {J}acobi-{P}erron type algorithms, and {D}irichlet groups},
  journal = {Acta Arith.},
  volume = {203},
  number = {1},
  pages = {27--48},
  year = {2022}
}

@incollection{Lean4,
  author = {de Moura, L. and Ullrich, S.},
  title = {The {L}ean 4 theorem prover and programming language},
  booktitle = {Automated Deduction -- {CADE} 28},
  series = {Lecture Notes in Comput. Sci.},
  volume = {12699},
  pages = {625--635},
  publisher = {Springer, Cham},
  year = {2021}
}

@inproceedings{Mathlib2020,
  author = {{The mathlib Community}},
  title = {The {L}ean mathematical library},
  booktitle = {Proceedings of the 9th ACM SIGPLAN International Conference on Certified Programs and Proofs (CPP 2020)},
  pages = {367--381},
  publisher = {ACM, New York},
  year = {2020}
}

\vspace{.5cm}

\end{document}